\documentclass[12pt]{article}
\usepackage{graphicx}
\usepackage{amsmath}
\usepackage{amsfonts}
\usepackage{amssymb}

\newtheorem{theorem}{Theorem}

\newtheorem{corollary}[theorem]{Corollary}

\newtheorem{lemma}[theorem]{Lemma}

\newtheorem{proposition}[theorem]{Proposition}

\newenvironment{proof}[1][Proof]{\textbf{#1.} }{\ \rule{0.5em}{0.5em}}

\begin{document}

\title{Genus one 1-bridge knots and Dunwoody
manifolds\footnote{Work performed under the auspices of
G.N.S.A.G.A. of C.N.R. of Italy and supported by the University of
Bologna, funds for selected research topics.}}
\author{Luigi Grasselli \and Michele Mulazzani}

\maketitle
\begin{abstract}
{In this paper we show that all 3-manifolds of a family introduced
by M. J. Dunwoody are cyclic coverings of lens spaces (eventually
$\bf S^3$), branched over genus one 1-bridge knots. As a
consequence, we give a positive answer to the Dunwoody conjecture
that all the elements of a wide subclass are cyclic coverings of
$\bf S^3$ branched over a knot. Moreover, we show that all
branched cyclic coverings of a 2-bridge knot belong to this
subclass; this implies that the fundamental group of each branched
cyclic covering of a 2-bridge knot admits a geometric cyclic
presentation.
}

\end{abstract}

\section{Introduction and preliminaries}

The problem of determining if a balanced presentation of a group
is geometric (i.e. induced by a Heegaard diagram of a closed
orientable 3-manifold) is quite important within geometric
topology and has been deeply investigated by many authors (see
\cite{Gr}, \cite{Mo}, \cite{Ne}, \cite{OS1}, \cite{OS2},
\cite{OS3}, \cite{St}); further, the connections between branched
cyclic coverings of links and cyclic presentations of groups
induced by suitable Heegaard diagrams have been recently pointed
out in several papers (see \cite{BKM}, \cite{CHK},
\cite{CHR}, \cite{Du}, \cite{HKM1}, \cite{HKM2},
\cite{Ki}, \cite{KV}, \cite{KKV}, \cite{MR},
\cite{VK}). In order to investigate these connections, M.J.
Dunwoody introduces in \cite{Du} a class of planar, 3-regular
graphs endowed with a cyclic symmetry. Each graph is defined by a
6-tuple of integers; if this 6-tuple satisfies suitable conditions
(admissible 6-tuple), the graph uniquely defines a Heegaard
diagram such that the presentation of the fundamental group of the
represented manifold is cyclic. This construction gives rise to a
wide class of closed orientable 3-manifolds (Dunwoody manifolds),
depending on 6-tuples of integers and admitting geometric cyclic
presentations for their fundamental groups. Our main result is
that each Dunwoody manifold is a cyclic covering of a lens space
(eventually the 3-sphere), branched over a genus one 1-bridge
knot. As a direct consequence, the Dunwoody manifolds belonging to
a wide subclass are proved to be cyclic coverings of $\bf S^3$,
branched over suitable knots, thus giving a positive answer to a
conjecture of Dunwoody \cite{Du}. Moreover, we show that all
branched cyclic coverings of knots with classical (i.e. genus
zero) bridge number two belong to this subclass; as a corollary,
the fundamental group of each branched cyclic covering of a
2-bridge knot admits a geometric cyclic presentation.

For the theory of Heegaard splittings of 3-manifolds, and in
particular for Singer moves on Heegaard diagrams realizing the
homeomorphism of the represented manifolds, we refer to \cite{He}
and \cite{Si}. For the theory of cyclically presented groups, we
refer to \cite{Jo}.

We recall that a finite balanced presentation of a group
$<x_1,\ldots,x_n\vert r_1,\ldots,r_n>$ is said to be a {\it cyclic
presentation\/} if there exists a word $w$ in the free group $F_n$
generated by $x_1,\ldots,x_n$ such that the relators of the
presentation are $r_k=\theta_n^{k-1}(w)$, $k=1,\ldots,n$, where
$\theta_n :F_n\to F_n$ denotes the automorphism defined by
$\theta_n (x_i)=x_{i+1}$ (mod $n$), $i=1,\ldots,n$. Let us denote
this cyclic presentation (and the related group) by the symbol
$G_n(w)$, so that: $$G_n(w)=<x_1,x_2,\ldots,x_n\vert
w,\theta_n(w),\ldots,\theta_n^{n-1}(w)>.$$

A group is said to be cyclically presented if it admits a cyclic
presentation. We recall that the {\it exponent-sum\/} of a word
$w\in F_n$ is the integer $\varepsilon_w$ given by the sum of the
exponents of its letters; in other terms,
$\varepsilon_w=\upsilon(w)$ where $\upsilon:F_n\to \bf Z$ is the
homomorphism defined by $\upsilon(x_i)=1$ for each $1\le i\le n$.

Following \cite{Ha}, we recall the definition of genus $g$ bridge
number of a link, which is a generalization of the classical
concept of bridge number for links in $\bf S^3$ (see \cite{Do}).

A set of mutually disjoint arcs $\{t_1,\ldots,t_n\}$ properly
embedded in a handlebody $U$ is trivial if there is a set of
mutually disjoint discs $D=\{D_1,\ldots,D_n\}$ such that $t_i\cap
D_i=t_i\cap\partial D_i=t_i$, $t_i\cap D_j=\emptyset$ and
$\partial D_i-t_i\subset\partial U$ for $1\le i,j\le n$ and $i\neq
j$. Let $U_1$ and $U_2$ be the two handlebodies of a Heegaard
splitting of the closed orientable 3-manifold $M$ and let $T$ be
their common surface: a link $L$ in $M$ is in $n$-bridge position
with respect to $T$ if $L$ intersects $T$ transversally and if the
set of arcs $L\cap U_i$ has $n$ components and is trivial both in
$U_1$ and in $U_2$. A link in 1-bridge position is obviously a
knot.

The genus $g$ bridge number of a link $L$ in $M$, $b_g(L)$,  is
the smallest integer $n$ for which $L$ is in $n$-bridge position
with respect to some genus $g$ Heegaard surface in $M$. If the
genus $g$ bridge number of a link $L$ is $b$, we say that $L$ is a
genus $g$ $b$-bridge link or simply a ($g,b$)-link. Of course, the
genus $g$ bridge number of a link in a manifold of Heegaard genus
$g'$ is defined only for $g\ge g'$ and the genus 0 bridge number
of a link in $\bf S^3$ is the classical bridge number. Moreover, a
($g,1$)-link is a knot, for each $g\ge 0$.

In what follows, we shall deal with ($1,1$)-knots, i.e. knots in
$\bf S^3$ or in lens spaces. This class of knots is very important
in the light of some results and conjectures involving Dehn
surgery on knots (see \cite{Be}, \cite{Ga1}, \cite{Ga2},
\cite{W1}, \cite{W2}, \cite{W3}). Notice that the class of
($1,1$)-knots in $\bf S^3$ contains all torus knots (trivially)
and all 2-bridge knots (i.e. $(0,2)$-knots) \cite{MS}.

\section{Dunwoody manifolds}

Let us sketch now the construction of Dunwoody manifolds given in
\cite{Du}. Let $a,b,c,n$ be integers such that $n>0$, $a,b,c\ge 0$
and $a+b+c>0$. Let $\Gamma=\Gamma(a,b,c,n)$ be the planar regular
trivalent graph drawn in Figure 1.

\bigskip

\begin{figure}[bht]
 \begin{center}
 \includegraphics*[totalheight=7cm]{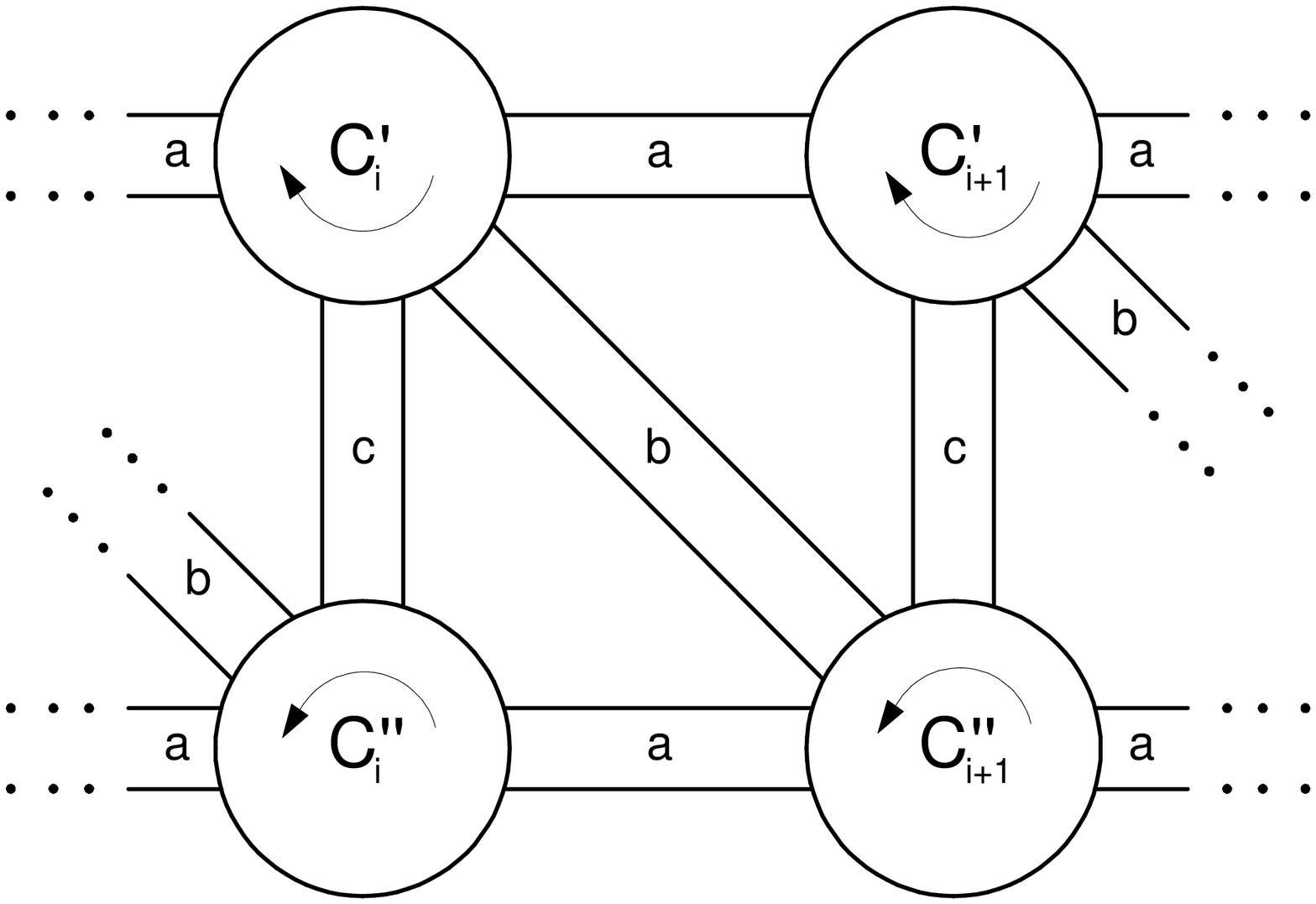}
 \end{center}
 \caption{The graph $\Gamma$.}

 \label{Fig. 1}

\end{figure}

It contains $n$ upper cycles
$C'_1,\ldots,C'_n$ and $n$ lower cycles $C''_1,\ldots,C''_n$, each
having $d=2a+b+c$ vertices. For each $i=1,\ldots,n$, the cycle
$C'_i$ (resp. $C''_i$) is connected to the cycle $C'_{i+1}$ (resp.
$C''_{i+1}$) by $a$ parallel arcs, to the cycle $C''_{i}$ by $c$
parallel arcs and to the cycle $C''_{i+1}$ by $b$ parallel arcs
(assume $n+1=1$). We set ${\cal C'}=\{C'_1,\ldots,C'_n\}$ and
${\cal C''}=\{C''_1,\ldots,C''_n\}$. Moreover, denote by $A'$
(resp. $A''$) the set of the arcs of $\Gamma$ belonging to a cycle
of $\cal C'$ (resp. $\cal C''$) and by $A$ the set of the other
arcs of the graph. The one-point compactification of the plane
leads to a 2-cell embedding of the graph $\Gamma$ in $\bf S^2$; it
is evident that the graph is invariant with respect to a rotation
$\rho_n$ of the sphere by $2\pi /n$ radians along a suitable axis
intersecting $\bf S^2$ in two points not belonging to the graph.
Obviously, $\rho_n$ sends $C'_i$ to $C'_{i+1}$ and $C''_i$ to
$C''_{i+1}$ (mod $n$), for each $i=1,\ldots,n$.

By cutting the sphere along all $C'_i$ and $C''_i$ and by removing
the interior of the corresponding discs, we obtain a sphere with
$2n$ holes. Let now $r$ and $s$ be two new integers; give a
clockwise (resp. counterclockwise) orientation to the cycles of
$\cal C'$ (resp. of $\cal C''$) and label their vertices from $1$
to $d$, in accordance with these orientations (see Figure 2) so
that:
\begin{itemize}
\item [-] the vertex 1 of each $C'_i$ is the endpoint of the first arc of
$A$ connecting $C'_i$ with $C'_{i+1}$;
\item [-] the vertex $1-r$ (mod $d$) of each $C''_i$ is the endpoint of the
first arc of $A$ connecting $C''_i$ with $C''_{i+1}$.
\end{itemize}
Then glue the cycle $C'_i$ with the cycle $C''_{i-s}$ (mod $n$) so
that equally labelled vertices are identified together.

\bigskip

\begin{figure}[bht]
 \begin{center}
 \includegraphics*[totalheight=8cm]{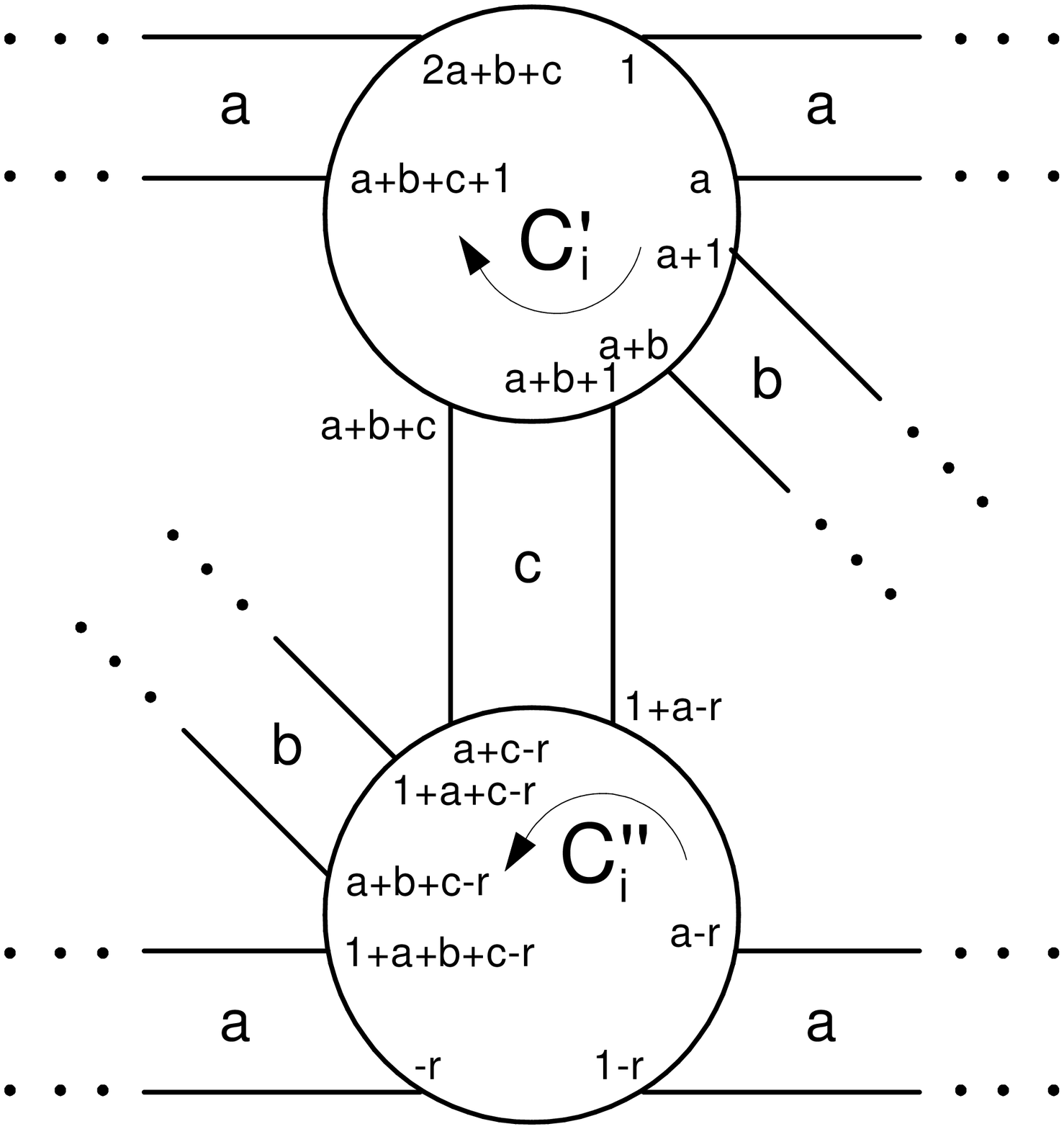}
 \end{center}
 \caption{}

 \label{Fig. 2}

\end{figure}

It is evident by construction that the integers $r$ and $s$ can be
taken mod $d$ and mod $n$ respectively. Denote by $\cal S$ the set
of all the 6-tuples $(a,b,c,n,r,s)\in \bf Z^6$ such that $n>0$,
$a,b,c\ge 0$ and $a+b+c>0$.

The described gluing gives rise to an orientable surface $T_n$ of
genus $n$ and the $nd$ arcs belonging to $A$ are pairwise
connected through their endpoints, realizing $m$ cycles
$D_1,\ldots,D_m$ on $T_n$. It is straightforward that the cut of
$T_n$ along the $n$ cycles $C_i=C'_i=C''_{i-s}$ does not
disconnect the surface. Set ${\cal C}=\{C_1,\ldots,C_n\}$ and
${\cal D}=\{D_1,\ldots,D_m\}$.

If $m=n$ and if the cut along the cycles of ${\cal D}$ does not
disconnect $T_n$, then the two systems of meridian curves $\cal C$
and $\cal D$ in $T_n$ represent a genus $n$ Heegaard diagram of a
closed orientable 3-manifold, which is completely determined by
the 6-tuple. Each manifold arising in this way is called a {\it
Dunwoody manifold\/}.

Thus, we define to be {\it admissible\/} the 6-tuples
$(a,b,c,n,r,s)$ of $\cal S$ satisfying the following conditions:

\begin{itemize}
\item [(1)] the set $\cal D$ contains exactly $n$ cycles;
\item [(2)] the surface $T_n$ is not disconnected by the cut
along the cycles of $\cal D$.
\end{itemize}

The ``open'' Heegaard diagram $\Gamma$ and the Dunwoody manifold
associated to the admissible 6-tuple $\sigma$ will be denoted by
$H(\sigma)$ and $M(\sigma)$ respectively.

\medskip

\noindent{\bf Remark 1.} It is easy to see that not all the
6-tuples in $\cal S$ are admissible. For example, the 6-tuples
$(a,0,a,1,a,0)$, with $a\ge 1$, give rise to exactly $a$ cycles in
$\cal D$; thus, they are not admissible if $a>1$. The 6-tuples
$(1,0,c,1,2,0)$ are not admissible if $c$ is even, since, in this
case, we obtain exactly one cycle $D_1$, but the cut along it
disconnects the torus $T_1$.

\medskip

Consider now a 6-tuple $\sigma\in\cal S$. The graph $\Gamma$
becomes, via the gluing quotient map, a regular 4-valent graph
denoted by $\Gamma'$ embedded in $T_n$. Its vertices are the
intersection points of the spaces $\Omega=\cup_{i=1}^n C_i$ and
$\Lambda=\cup_{j=1}^m D_j$; hence they inherit the labelling of
the corresponding glued vertices of $\Gamma$. Since the gluing of
the cycles of $\cal C'$ and $\cal C''$ is invariant with respect
to the rotation $\rho_n$, the group ${\cal G}_n=<\rho_n>$
naturally induces a cyclic action of order $n$ on $T_n$ such that
the quotient $T_1=T_n/{\cal G}_n$ is homeomorphic to a torus. The
labelling of the vertices of $\Gamma'$ is invariant under the
rotation $\rho_n$ and $\rho_n(C_i)=C_{i+1}$ (mod $n$). We are
going to show that, if the 6-tuple is admissible, this last
property also holds for the cycles of $\cal D$.


\begin{lemma} \label{Lemma 0'} a) Let $\sigma=(a,b,c,n,r,s)$ be an admissible 6-tuple.
Then $\rho_n$ induces a cyclic permutation on the curves of $\cal
D$. Thus, if $D$ is a cycle of $\cal D$, then ${\cal
D}=\{\rho_n^{k-1}(D)\vert k=1,\ldots,n\}$.

b) If $(a,b,c,n,r,s)$ is admissible, then also $(a,b,c,1,r,0)$ is
admissible and the Heegaard diagram $H(a,b,c,1,r,0)$ is the
quotient of the Heegaard diagram $H(a,b,c,n,r,s)$ respect to
${\cal G}_n$.
\end{lemma}


\begin{proof} a) First of all, note that $\rho_n(\Lambda)=\Lambda$; thus
the group ${\cal G}_n$ also acts on the spaces $T_n-\Lambda$ and
$\Lambda$ (and hence on the set $\cal D$). If the 6-tuple $\sigma$
is admissible, then $T_n-\Lambda$ is connected, and hence the
quotient $(T_n-\Lambda)/{\cal G}_n=T_n/{\cal G}_n-\Lambda/{\cal
G}_n$ must be connected too. This implies that $\Lambda/{\cal
G}_n$ has a unique connected component. Since $\Lambda$ has
exactly $n$ connected components, the cyclic group ${\cal G}_n$ of
order $n$ defines a simply transitive cyclic action on the cycles
of $\cal D$.

b) Let $C,D\subset T_1$ the two curves $C=\Omega/{\cal G}_n$ and
$D=\Lambda/{\cal G}_n$. Then, the two systems of curves ${\cal
C}=\{C\}$ and ${\cal D}=\{D\}$ on $T_1$ define a Heegaard diagram
of genus one. The graph $\Gamma_1$ corresponding to
$\sigma_1=(a,b,c,1,r,0)$ is the quotient of the graph $\Gamma_n$
corresponding to $\sigma=(a,b,c,n,r,s)$, respect to ${\cal G}_n$.
Moreover, the gluings on $\Gamma_n$ are invariant respect to
$\rho_n$. Therefore, the gluings on $\Gamma_1$ give rise to the
Heegaard diagram above. This show that the 6-tuple $\sigma_1$ is
admissible and obviously $H(a,b,c,1,r,0)$ is the quotient of
$H(a,b,c,n,r,s)$ respect to ${\cal G}_n$.
\end{proof}

\medskip

\noindent{\bf Remark 2.} More generally, given two positive
integer $n$ and $n'$ such that $n'$ divides $n$, if
$(a,b,c,r,n,s)$ is admissible, then $(a,b,c,r,n',b)$ is admissible
too. Moreover, the Heegaard diagram $H(a,b,c,r,n',b)$ is the
quotient of $H(a,b,c,r,n,b)$ respect to the action of a cyclic
group of order $n/n'$.

\medskip

It is easy to see that, for admissible 6-tuples, each cycle in
$\cal D$ contains $d$ vertices with different labels and is
composed by exactly $d$ arcs of $\Gamma$ (in fact, $2a$ horizontal
arcs, $b$ oblique arcs and $c$ vertical arcs).

An important consequence of point a) of Lemma \ref{Lemma 0'} is
that, if $\sigma$ is an admissible 6-tuple, the presentation of
the fundamental group of $M(\sigma)$ induced by the Heegaard
diagram $H(\sigma)$ is cyclic.

To see this, let $v$ be the vertex belonging to the cycle $C_1$
and labelled by $a+b+1$; denote by $D_1$ the curve of $\cal D$
containing $v$ and by $v'$ the vertex of $C'_1$ corresponding to
$v$. Orient the arc $e'\in A$ of the graph $\Gamma$ containing
$v'$ so that $v'$ is its first endpoint and orient the curve $D_1$
in accordance with the orientation of this arc. Now, set
$D_k=\rho_n^{k-1}(D_1)$, for each $k=1,\ldots,n$; the orientation
on $D_1$ induces, via $\rho_n$, an orientation also on these
curves. Moreover, these orientation on the cycles of $\cal D$
induce an orientation on the arcs of the graph $\Gamma$ belonging
to $A$. By orienting the arcs of $C'$ and $C''$ in accordance with
the fixed orientations of the cycles $C'_i$ and $C''_i$, the graph
$\Gamma$ becomes an oriented graph, whose orientation is invariant
under the action of the group ${\cal G}_n$. Let us define to be
{\it canonical\/} this orientation of $\Gamma$.

Let now $w\in F_n$ be the word obtained by reading the oriented
arcs $e_1=e',e_2,\ldots,e_d$ of $\Gamma$ corresponding to the
oriented cycle $D_1$, starting from the vertex $v'$. The letters
of $w$ are in one-to-one correspondence with the oriented arcs
$e_h$; more precisely, the letter of $w$ corresponding to $e_h$ is
$x_i$ if $e_h$ comes out from the cycle $C'_i$ and is $x_i^{-1}$
if $e_h$ comes out from the cycle $C''_{i-s}$. Note that the word
$\theta_n^{k-1}(w)$ in the cyclic presentation $G_n(w)$ is
obtained by reading the cycle $D_k$ along the given orientation,
for $1\le k\le n$ (roughly speaking, the automorphism $\theta_n$
is ``geometrically'' realized by $\rho_n$).

This proves that each admissible 6-tuple $\sigma$ uniquely
defines, via the associated Heegaard diagram $H(\sigma)$, a word
$w=w(\sigma)$ and a cyclic presentation $G_n(w)$ for the
fundamental group of the Dunwoody manifold $M(\sigma)$. Note that
the sequence of the exponents in the word $w(\sigma)$, and hence
its exponent-sum $\varepsilon_{w(\sigma)}$, only depends on the
integers $a,b,c,r$.

Let us consider now the Dunwoody manifolds $M(a,b,c,n,r,s)$ with
$n=1$ (and hence $s=0$), which arises from a genus one Heegaard
diagram.


\begin{proposition} \label{Proposition 1} Let $(a,b,c,1,r,0)$ be
an admissible 6-tuple and let $w=w(a,b,c,1,r,0)$ be the associated
word. Then the Dunwoody manifold $M(a,b,c,1,r,0)$ is homeomorphic
to:

i) $\bf S^3$, if $\varepsilon_w=\pm 1$;

ii) $\bf S^1\times S^2$, if $\varepsilon_w=0$;

iii) a lens space $L(\alpha,\beta)$ with $\alpha=\vert
\varepsilon_w\vert$, if $\vert\varepsilon_w\vert>1$.
\end{proposition}


\begin{proof} From $n=1$ we obtain $w\in F_1\cong {\bf Z}\cong
<x\vert\emptyset>$. Thus, $\pi_1(M)\cong G_1(w)\cong <x\vert
x^{\varepsilon_w}>\cong {\bf Z}_{\vert\varepsilon_w\vert}$.
\end{proof}

\medskip

\noindent{\bf Example 1.} The Dunwoody manifolds $M(0,0,1,1,0,0)$,
$M(1,0,0,1,1,0)$ and $M(0,0,c,1,r,0)$, with $c,r$ coprime, are
homeomorphic to $\bf S^3$, $\bf S^1\times S^2$ and to the lens
space $L(c,r)$, respectively. Moreover, all lens spaces also arise
with $a\ne 0$ ; in fact, for each $a>0$, $M(a,0,c,1,a,0)$ is
homeomorphic with the lens space $L(c,a)$, if $a$ and $c$ are
coprime, since it is easy to see that $H(a,0,c,1,a,0)$ can be
transformed into the canonical genus one Heegaard diagram of
$L(c,a)$ by Singer moves of type IB.

\medskip

Let us see now how the admissibility conditions for the 6-tuples
of $\cal S$ can be given in terms of labelling of the vertices of
$\Gamma'$, belonging to the curve $D_1\in\cal D$. With this aim,
consider the following properties for a 6-tuple $\sigma\in\cal S$:
\begin{itemize}
\item [(i')] the set of the labels of the vertices belonging to
the cycle $D_1$ is the set of all integers from $1$ to $d$;
\item [(ii')] the vertices of the cycle $D_1$ have different
labels.
\end{itemize}

It is easy to see that, if a 6-tuple $\sigma\in\cal S$ is
admissible, then it satisfies (i') and (ii'). On the other side,
if a 6-tuple $\sigma\in\cal S$ satisfies (i') and (ii'), then the
curves $\rho_n^{k-1}(D_1)\in\cal D$, with $k=1,\ldots,n$, which
are all different from each other, are precisely the curves of $\cal
D$. Thus, $\cal D$ has exactly $n$ curves and they are cyclically
permutated by $\rho_n$. However, this does not imply that $\sigma$
is admissible; for example, the 6-tuple $(1,0,2,1,2,0)$ satisfies
(i') and (ii'), but it is not admissible (see Remark 1). Note
that, for $n=1$, property (ii') always holds, while condition (i')
holds if and only if $\cal D$ has a unique cycle.

If a 6-tuple satisfies property (i'), then ${\cal G}_n$ acts
transitively (not necessarily simply) on $\cal D$, and hence it is
possible to induce an orientation (which is still said to be
canonical) on the cycles of $\cal D$ and on the graph $\Gamma$, by
extending, via $\rho_n$, the orientation of $D_1$ to the other
cycles of $\cal D$.

Property (i') implies that the cycles of $\cal D$ naturally induce
a cyclic permutation on the set ${\cal N}=\{1,\ldots,d\}$ of the
vertex labels. In fact, by walking along these canonically
oriented cycles, starting from an arbitrary vertex $\bar v$
labelled $j$, one sequentially meets $d$ vertices (whose labels
are different from each other), and then a new vertex $\bar v'$
labelled $j$ which can be different from $\bar v$. The sequence of
the labellings of these $d$ consecutive vertices defines the
cyclic permutation on $\cal N$. Further, each cycle of $\cal D$
precisely contains $d'=ld$ arcs, with $l\ge 1$, and $l=1$ if and
only if the 6-tuple satisfies (ii') too. Moreover, property (i')
is independent from the integers $n$ and $s$; hence, given two
6-tuples $\sigma=(a,b,c,n,r,s)$ and $\sigma'=(a,b,c,n',r,s)$, then
$\sigma$ satisfies (i') if and only if $\sigma'$ satisfies (i').

Let now $\sigma$ be a 6-tuple satisfying (i') and suppose that
$\Gamma$ is canonically oriented. An arc of $\Gamma$ belonging to
$A$ is said to be of type I if it is oriented from a cycle of
$\cal C'$ to a cycle of $\cal C''$, of type II if it is oriented
from a cycle of $\cal C''$ to a cycle of $\cal C'$ and of type III
otherwise (it joins cycles of $\cal C'$ or cycles of $\cal C''$).
Moreover, the arc is said to be of type I' if it is oriented from
a cycle $C'_i$ (resp. $C''_i$) to a cycle $C'_{i+1}$ (resp.
$C''_{i+1}$), of type II' if it is oriented from a cycle
$C'_{i+1}$ (resp. $C''_{i+1}$) to a cycle $C'_i$ (resp. $C''_i$)
and of type III' otherwise (it joins $C'_i$ with $C''_i$). Let
$\Delta$ be the set of the first $d$ arcs of $D_1$, following the
canonical orientation, starting from the arc coming out from the
vertex $v'$ of $C'_1$ labelled $a+b+1$. Obviously, the set
$\Delta$ contains all the arcs of $D_1$ if and only if the 6-tuple
$\sigma$ also satisfies (ii').

Now, denote by $p'_{\sigma}$ (resp. $p''_{\sigma}$) the number of
the arcs of type I (resp. of type II) of $\Delta$ and set
$p_{\sigma}=p'_{\sigma}-p''_{\sigma}$. Similarly, denote by
$q'_{\sigma}$ (resp. $q''_{\sigma}$) the number of the arcs of
type I' (resp. of type II') of $\Delta$ and set
$q_{\sigma}=q'_{\sigma}-q''_{\sigma}$. Note that $p_{\sigma}$ has
the same parity of $b+c$ and $q_{\sigma}$ has the same parity of
$2a+b$ and hence of $b$. It is evident that $p_{\sigma}$ and
$q_{\sigma}$ only depend on the integers $a,b,c,r$.

The integers $p_{\sigma}$ and $q_{\sigma}$ give an useful tool for
verifying condition (ii'). In fact, suppose to walk along the
canonically oriented cycle $D_j$ of $\cal D$, starting from a
vertex $\bar v$ and let $C_i$ be the cycle of $\cal C$ containing
$\bar v$. If $\bar v'$ is the first vertex with the same label of
$\bar v$ and if $C_{i'}$ is the cycle of $\cal C$ containing $\bar
v'$, we have $i'=i+q_{\sigma}+sp_{\sigma}$. Thus, the cycle $D_j$
contains $d$ arcs if and only if $q_{\sigma}+sp_{\sigma}\equiv 0$
(mod $n$). This proves that the 6-tuple satisfies (ii'). Thus,
(i') and (ii') are respectively, in a different language,
conditions (i) and (ii) of Theorem 2 of \cite{Du}, which gives a
necessary and sufficient condition for a 6-tuple to be admissible
when $d$ is odd. In fact, we have the following result:

\begin{lemma} \label{Lemma 2'} (\cite{Du}, Theorem 2) Let
$\sigma=(a,b,c,n,r,s)$ be a 6-tuple with $d=2a+b+c$ odd. Then
$\sigma$ is admissible if and only if it satisfies (i') and (ii').
\end{lemma}

\noindent{\bf Remark 3.} This result does not hold when $d$ is
even. In fact, the 6-tuples $(1,0,c,1,2,0)$, with $c$ even,
satisfy (i') and (ii'), but they are not admissible, as pointed
out in Remark 1.

\medskip

An immediate consequence of Lemma \ref{Lemma 2'} is the following result:

\begin{corollary} \label{Corollary 2''} Let $\sigma=(a,b,c,n,r,s)$ be a 6-tuple with
$d=2a+b+c$ odd and $n=1$. Then $\sigma$ is admissible if and only
if $\cal D$ has a unique cycle.
\end{corollary}

\begin{proof}If $\sigma$ is admissible, then it is straightforward that
$\cal D$ has a unique cycle. Vice versa, if $\cal D$ has a unique
cycle, then (i') holds. Since $n=1$ implies (ii'), the result is a
direct consequence of the above lemma.
\end{proof}

The parameter $p_{\sigma}$ associated to an admissible 6-tuple
$\sigma$ is strictly related to the word $w(\sigma)$ associated to
$\sigma$. In fact, we have:


\begin{lemma} \label{Lemma 2} Let $\sigma=(a,b,c,n,r,s)$ be an admissible 6-tuple,
$w=w(\sigma)$ the associated word and $\varepsilon_w$ its
exponent-sum. Then $$p_{\sigma}=\varepsilon_w.$$
\end{lemma}

\begin{proof} Since $\sigma$ is admissible, the arcs of $\Delta$ are
precisely the arcs of $D_1$. Let $e_1,e_2,\ldots,e_d$ be the
sequence of these arcs, following the canonical orientation on
$D_1$, and let $w=\prod_{h=1}^d x_{i_h}^{u_h}$, with
$u_h\in\{+1,-1\}$. We have: $\varepsilon_w=\sum_{h=1}^d
u_h=1/2\sum_{h=1}^d(u_h+u_{h+1})$,where $d+1=1$. Since
$u_h+u_{h+1}=+2$ if $e_h$ is of type I, $u_h+u_{h+1}=-2$ if $e_h$
is of type II and $u_h+u_{h+1}=0$ if $e_h$ is of type III, the
result immediately follows.
\end{proof}

In \cite{Du} Dunwoody investigates a wide subclass of manifolds
$M(\sigma)$ such that $p_{\sigma}=\pm 1$ and he conjectures that
all the elements of this subclass are cyclic coverings of $\bf
S^3$ branched over knots. In the next chapter this conjecture will
be proved as a corollary of a more general theorem.

\section{Main results}

The following theorem is the main result of this paper and shows
how the cyclic action on the Heegaard diagrams naturally extends
to a cyclic action on the associated Dunwoody manifolds, which
turn out to be cyclic coverings of $\bf S^3$ or of lens spaces,
branched over suitable knots.


\begin{theorem} \label{Theorem 3} Let $\sigma=(a,b,c,n,r,s)$ be an admissible 6-tuple,
with $n>1$. Then the Dunwoody manifold $M=M(a,b,c,n,r,s)$ is the
$n$-fold cyclic covering of the manifold $M'=M(a,b,c,1,r,0)$,
branched over a genus one 1-bridge knot $K=K(a,b,c,r)$ only
depending on the integers $a,b,c,r$. Further, $M'$ is homeomorphic
to:

i) $\bf S^3$, if $p_{\sigma}=\pm 1$,

ii) $\bf S^1\times S^2$, if $p_{\sigma}=0$,

iii) a lens space $L(\alpha,\beta)$ with $\alpha=\vert
p_{\sigma}\vert$, if $\vert p_{\sigma}\vert>1$.
\end{theorem}

\begin{proof} Since the two systems of curves ${\cal
C}=\{C_1,\ldots,C_n\}$ and ${\cal D}=\{D_1,\ldots,D_n\}$ on $T_n$
define a Heegaard diagram of $M$, there exist two handlebodies
$U_n$ and $U'_n$ of genus $n$, with $\partial U_n=\partial
U'_n=T_n$, such that $M=U_n\cup U'_n$. Let now ${\cal G}_n$ be the
cyclic group of order $n$ generated by the homeomorphism $\rho_n$
on $T_n$. The action of ${\cal G}_n$ on $T_n$ extends to both the
handlebodies $U_n$ and $U'_n$ (see \cite{RZ}), and hence to the
3-manifold $M$. Let $B_1$ (resp. $B'_1$) be a disc properly
embedded in $U_n$ (resp. in $U'_n$) such that $\partial B_1=C_1$
(resp. $\partial B'_1=D_1$). Since $\rho_n(C_i)=C_{i+1}$ and
$\rho_n(D_i)=D_{i+1}$ (mod $n$), the discs $B_k=\rho_n^{k-1}(B_1)$
(resp. $B'_k=\rho_n^{k-1}(B'_1)$), for $k=1,\ldots,n$, form a
system of meridian discs for the handlebody $U_n$ (resp. $U'_n$).
By arguments contained in \cite{Zi}, the quotients $U_1=U_n/{\cal
G}_n$ and $U'_1=U'_n/{\cal G}_n$ are both handlebody orbifolds
topologically homeomorphic to a genus one handlebody with one arc
trivially embedded as its singular set with a cyclic isotropy
group of order $n$. The intersection of these orbifolds is a
2-orbifold with two singular points of order $n$, which is
topologically the torus $T_1=T_n/{\cal G}_n$; the curve $C$ (resp.
$D$), which is the image via the quotient map of the curves $C_i$
(resp. of the curves $D_i$), is non-homotopically trivial in
$T_1$. These curves, each of which is a fundamental system of
curves in $T_1$, define a Heegaard diagram of $M'$ (induced by
$H(a,b,c,1,r,0)$). The union of the orbifolds $U_1$ and $U'_1$ is
a 3-orbifold topologically homeomorphic to $M'$, having a genus
one 1-bridge knot $K\subset M'$ as singular set of order $n$.
Thus, $M'$ is homeomorphic to $M/{\cal G}_n$ and hence $M$ is the
$n$-fold cyclic covering  of $M'$, branched over $K$. Since the
handlebody orbifolds and their gluing only depend on $a,b,c,r$,
the same holds for the branching set $K$. The homeomorphism type
of $M'$ follows from Proposition \ref{Proposition 1} and Lemma
\ref{Lemma 2}.
\end{proof}

\medskip

\noindent{\bf Remark 4.} More generally, given two positive
integers $n$ and $n'$ such that $n'$ divides $n$, if
$(a,b,c,r,n,s)$ is admissible, then the Dunwoody manifold
$M(a,b,c,n,r,s)$ is the $n/n'$-fold cyclic covering of the
manifold $M'=M(a,b,c,n',r,s)$, branched over an $(n',1)$-knot in
$M'$.

\medskip

\noindent{\bf Example 2.} The Dunwoody manifolds $M(0,0,1,n,0,0)$,
$M(1,0,0,n,1,0)$ and $M(0,0,c,n,r,0)$, with $c,r$ coprime, are
$n$-fold cyclic coverings of the manifolds $\bf S^3$, $\bf
S^1\times S^2$ and $L(c,r)$ respectively, branched over a trivial
knot. In fact, these Dunwoody manifolds are the connected sum of
$n$ copies of $\bf S^3$, $\bf S^1\times S^2$ and $L(c,r)$
respectively.

\medskip

Let us consider now the class of the Dunwoody manifolds
$M_n=M(a,b,c,n,r,s)$ with $p=\pm 1$ (and hence $d$ odd) and
$s=-pq$. Many examples of these manifolds appear in Table 1 of
\cite{Du}, where it was conjectured that they are $n$-fold cyclic
coverings of $\bf S^3$, branched over suitable knots. The
following corollary of Theorem \ref{Theorem 3} proves this
conjecture.


\begin{corollary} \label{Corollary 4} Let $\sigma_1=(a,b,c,1,r,0)$
be an admissible 6-tuple with $p_{\sigma_1}=\pm 1$ and
$s=-p_{\sigma_1}q_{\sigma_1}$. Then the 6-tuple
$\sigma_n=(a,b,c,n,r,s)$ is admissible for each $n>1$ and the
Dunwoody manifold $M_n=M(a,b,c,n,r,s)$ is a $n$-fold cyclic
coverings of $\bf S^3$, branched over a genus one 1-bridge knot
$K\subset \bf S^3$, which is independent on $n$.
\end{corollary}

\begin{proof} Obviously $(a,b,c,1,r,s)=\sigma_1$. Since $\sigma_1$ is
admissible, it satisfies (i'). This proves that $\sigma_n$
satisfies (i'), for each $n>1$. Since
$s=-p_{\sigma_1}q_{\sigma_1}=-p_{\sigma_n}q_{\sigma_n}$ and
$p_{\sigma_n}=p_{\sigma_1}=\pm 1$, we obtain
$q_{\sigma_n}+sp_{\sigma_n}=0$, for each $n>1$, which implies
condition (ii) of Theorem 2 of \cite{Du}, or equivalently (ii').
Moreover, $d$ is odd, since
$[d]_2=[2a+b+c]_2=[b+c]_2=[p_{\sigma_n}]_2=[p_{\sigma_1}]_2=1$.
Thus, Lemma \ref{Lemma 2'} proves that $\sigma_n$ is admissible.
The final result is then a direct consequence of Theorem
\ref{Theorem 3}.
\end{proof}

\medskip

We point out that the above result has been independently obtained
by H. J. Song and S. H. Kim in \cite{SK}.

An interesting problem which naturally arises is that of
characterizing the set $\cal K$ of branching knots in $\bf S^3$
involved in Corollary \ref{Corollary 4}. The next theorem shows
that it contains all 2-bridge knots. We recall that a 2-bridge
knot is determined by two coprime integers $\alpha$ and $\beta$,
with $\alpha>0$ odd. The classification of 2-bridge knots and
links has been obtained by Schubert in \cite{Sc}. Since the
2-bridge knot of type $(\alpha,\beta)$ is equivalent to the
2-bridge knot of type $(\alpha,\alpha-\beta)$, then $\beta$ can be
assumed to be even.

\begin{theorem} \label{Theorem 5} The 6-tuple $\sigma_1=(a,0,1,1,r,0)$ with
$(2a+1,2r)=1$ is admissible. Moreover, if $s=-q_{\sigma_1}$, then
the 6-tuple $\sigma_n=(a,0,1,n,r,s)$ is admissible for each $n>1$
and the Dunwoody manifold $M_n=M(a,0,1,n,r,s)$ is the $n$-fold
cyclic covering of $\,\bf S^3$, branched over the 2-bridge knot of
type $(2a+1,2r)$. Thus, all branched cyclic coverings of 2-bridge
knots are Dunwoody manifolds.
\end{theorem}

\begin{proof} From $(2a+1,2r)=1$ it immediately follows that $\sigma_1$
has a unique cycle in $\cal D$. Since $d=2a+1$ is odd, Corollary
\ref{Corollary 2''} proves that $\sigma_1$ is admissible. Since
$p_{\sigma_n}=p_{\sigma_1}=+1$, all assumptions of Corollary
\ref{Corollary 4} hold; hence $\sigma_n$ is admissible for each
$n>1$ and $M_n$ is an $n$-fold cyclic covering of $\bf S^3$,
branched over a knot $K\subset \bf S^3$ which is independent on
$n$. In order to determine this knot, we can restrict our
attention to the case $n=2$. Note that
$[s]_2=[-q_{\sigma_1}]_2=[b]_2=0$ and hence $s$ is always even.
Thus, in the case $n=2$ we can suppose $s=0$. Let us consider now
the genus two Heegaard diagram $H(a,0,1,2,r,0)$. The sequence of
Singer moves \cite{Si} on this diagram, drawn in Figures 3--10 and
described in the Appendix of the paper, leads to the canonical
genus one Heegaard diagram of the lens space $L(2a+1,2r)$ (see
Figure 10). Since the representation of lens spaces (including
$\bf S^3$) as 2-fold branched coverings of $\bf S^3$ is unique
\cite{HR}, the result immediately holds.
\end{proof}

\medskip

\noindent{\bf Remark 5.} The Dunwoody manifold $M(a,0,1,n,r,s)$ of
Theorem \ref{Theorem 5} is homeomorphic to the Minkus manifold
$M_n(2a+1,2r)$ \cite{Mi} and the Lins-Mandel manifold
$S(n,2a+1,2r,1)$ \cite{LM, Mu}.

\medskip

An immediate consequence of Theorem \ref{Theorem 5} is:

\begin{corollary} \label{Corollary 6} The fundamental group of every branched cyclic
covering of a 2-bridge knot admits a cyclic presentation which is
geometric.
\end{corollary}

\noindent{\bf Remark 6.} In \cite{Mi} is shown that the
fundamental group of every branched cyclic covering of a 2-bridge
knot admits a cyclic presentation, but without pointing out that
this presentation is geometric.

\medskip

About the set $\cal K$ of knots in $\bf S^3$ involved in Corollary
\ref{Corollary 4}, we propose the following:

\medskip

\noindent{\bf Conjecture.} The set $\cal K$ contains all torus knots.

\medskip

If this conjecture is true, the set $\cal K$ contains knots with
an arbitrarily high number of bridges. Moreover, the conjecture
implies that every branched cyclic covering of a torus knot admits
a geometric cyclic presentation. The above conjecture is supported by
several cases contained in Table 1 of \cite{Du} (see \cite{SK}).
For example, the Dunwoody manifolds $M(1,2,3,n,4,4)$ (resp.
$M(1,3,4,n,5,5)$) are the $n$-fold branched cyclic coverings of
the 4-bridge torus knot $K(4,5)$ (resp. of the 5-bridge torus knot
$K(5,6)$).

\section{Appendix}

Now we show how to obtain, by means of Singer moves \cite{Si} on
the genus two Heegaard diagram $H(a,0,1,2,r,0)$ of Figure 3, the
canonical genus one Heegaard diagram of the lens space
$L(2a+1,2r)$ of Figure 10. The result will be achieved by a
sequence of exactly $a+4$ Singer moves: one of type ID, $a+2$ of
type IC and the final one of type III.

Figure 3 shows the open Heegaard diagram $H(a,0,1,2,r,0)$. Note
that, since $s=0$, the cycle $C'_1$ (resp. $C'_2$) is glued with
the cycle $C''_1$ (resp. $C''_2$). Let $D_1$ (resp. $D_2$) be the
cycle of the Heegaard diagram corresponding to the arc $e'$ (resp.
$e''$) coming out from the vertex $v'$ of $C'_1$ (resp. $v''$ of
$C'_2$) labelled $a+1$. Orient $D_1$ (resp. $D_2$) so that the arc
$e'$ (resp. $e''$) is oriented from up to down (resp. from down to
up). This orientation on $D_2$ is opposite to the canonical one
but, in this way, all the $2a$ arcs connecting $C'_1$ with $C'_2$
are oriented from $C'_1$ to $C'_2$ and all the $2a$ arcs
connecting $C''_1$ with $C''_2$ are oriented from $C''_2$ to
$C''_1$. The cycle $D_1$, besides the arc $e'$, has two arcs for
each $k=0,\ldots,a-1$, one joining the vertex of $C'_1$ labelled
$a+1-(1+2k)r$ with the vertex of $C'_2$ labelled $a+1+(1+2k)r$,
and the other one joining the vertex of $C''_2$ labelled
$a+1+(1+2k)r$ with the vertex of $C''_1$ labelled $a+1-(3+2k)r$.
The cycle $D_2$, besides the arc $a_2$, has two arcs for each
$k=0,\ldots,a-1$, one joining the vertex of $C'_1$ labelled
$a+1-(2+2k)r$ with the vertex of $C'_2$ labelled $a+1+(2+2k)r$,
the other joining the vertex of $C''_2$ labelled $a+1+2kr$ with
the vertex of $C''_1$ labelled $a+1-(2+2k)r$.

The first Singer move consists of replacing the curve $D_2$ with
the curve $D'_2=D_1+D_2$ (move of type ID of \cite{Si}) obtained
by isotopically approaching the arcs $e'$ and $e''$ until their
intersection becomes a small arc and by removing the interior of
this arc. The move is completed by shifting, with a small isotopy,
$D_1$ in $D'_1$ so that it becomes disjoint from $D'_2$.

The resulting Heegaard diagram is drawn in Figure 4. The new
$2a+1$ pairs of vertices obtained on $C'_1,C''_1,C'_2,C''_2$ are
labelled by simply adding a prime to the old label, while the
$4a+2$ pairs of fixed vertices keep their old labelling. Note that
each new vertex labelled $j'$ is placed, in the cycles
$C'_1,C''_1,C'_2$ and $ C''_2$, between the old vertices labelled
$j$ and $j+1$ respectively. The cycles $C'_2$ and $C''_2$ are no
longer connected by any arc, while the cycles $C'_1$ and $C''_1$
are connected by a unique arc (belonging to $D'_1$) joining the
vertex labelled $(a+1)'$ of $C'_1$ with the vertex labelled
$(a+1-r)'$ of $C''_1$. All the $3a$ arcs connecting $C'_1$ and
$C'_2$ are oriented from $C'_1$ to $C'_2$ and all the $3a$ arcs
which now connect $C''_1$ with $C''_2$ are oriented from $C''_2$
to $C''_1$. The cycle $D'_2$ contains exactly $4a+2$ arcs; more
precisely, for each $i=1,\ldots,2a+1$, it has one arc joining the
vertex labelled $i$ of $C'_1$ with the vertex labelled $2a+2-i$ of
$C'_2$ and one arc joining the vertex labelled $i$ of $C''_2$ with
the vertex labelled $2a+2-2r-i$ of $C'_2$. The cycle $D'_1$ is a
copy of the cycle $D_1$ and hence it contains $2a+1$ arcs. One of
these arcs connects $C'_1$ with $C''_1$; moreover, for each
$k=0,\ldots,a-1$, $D'_1$ has one arc joining the vertex of $C'_1$
labelled $(a+1-(1+2k)r)'$ with the vertex of $C'_2$ labelled
$(a+1+(1+2k)r)'$ and one arc joining the vertex of $C''_2$
labelled $(a+1+(1+2k)r)'$ with the vertex of $C''_1$ labelled
$(a+1-(3+2k)r)'$.

Now, apply to the diagram a Singer move of type IC, cutting along
the cycle $E$ (drawn in Figure 4) containing $C''_1$ and $C''_2$
and gluing the curve $C''_2$ of the resulting disc with $C'_2$.

The new Heegaard diagram obtained in this way shown in Figure 5.
It contains the new cycles $E'$ and $E''$, which are copies of the
cutting cycle $E$. These cycles replace $C'_2$ and $C''_2$ and
they both have a unique vertex ($w'$ and $w''$ respectively). The
cycle $E'$ (resp. $E''$) is connected with $C'_1$ (resp. with
$C''_1$) by an arc joining $w'$ (resp. $w''$) with the vertex
labelled $(a+1)'$ (resp. $(a+1-r)'$), oriented as in Figure 5. The
cycles $C'_1$ and $C''_1$ are joined by $3a+1$ arcs, all oriented
from $C'_1$ to $C''_1$; $2a+1$ of them belong to $D'_2$ and the
other $a$ belong to $D'_1$. More precisely, for each
$i=1,\ldots,2a+1$, there is an arc of $D'_2$ joining the vertex
labelled $i$ of $C'_1$ with the vertex labelled $i-2r$ of $C''_1$;
while, for each $k=0,\ldots,a-1$, there is an arc of $D'_1$
joining the vertex labelled $(a+1-(1+2k)r)'$ of $C'_1$ with the
vertex labelled $(a+1-(3+2k)r)'$ of $C''_1$.

Apply again a Singer move of type IC, cutting along the cycle
$F_1$ (drawn in Figure 5) containing $C''_1$ and $E''$ and gluing
the curve $C''_1$ of the resulting disc with $C'_1$.

The resulting Heegaard diagram is shown in Figure 6. It contains
the new cycles $F'_1$ and $F''_1$, which are copies of the cutting
cycle $F_1$. These cycles replace $C'_1$ and $C''_1$ and they both
have one vertex less. It is easy to see that the cycle $D'_2$ has
exactly the same $2a+1$ arcs connecting $F'_1$ and $F''_1$, all
oriented from $F'_1$ to $F''_1$; if the labelling of the vertices
of $F'_1$ and $F''_1$ is induced by the labelling of $F_1$ shown
in Figure 5, these arcs join pairs of vertices with the same
labelling of the previous step. The cycle $D'_1$ instead has one
arc less than in the previous step. In fact, it has $a-1$ arcs,
connecting $F'_1$ and $F''_1$, all oriented from $F'_1$ to $F''_1$
and joining the vertex labelled $(a+1-(1+2k)r)'$ of $F'_1$ with
the vertex labelled $(a+1-(3+2k)r)'$ of $F''_1$, for
$k=1,\ldots,a-1$.

Now, apply again a Singer move of type IC, cutting along the cycle
$F_2$ (drawn in Figure 6) containing $F''_1$ and $E''$ and gluing
the curve $F''_1$ of the resulting disc with $F'_1$.

The new Heegaard diagram only differs from the previous one for
containing one arc less in the cycle $D'_1$. By inductive
application of Singer moves of type IC, cutting along the cycle
$F_h$ (drawn in Figure 7) containing $F''_{h-1}$ and $E''$ and
gluing the curve $F''_{h-1}$ of the resulting disc with
$F'_{h-1}$, we obtain, for $h=a$, the situation shown in Figure 8,
where the cycle $D'_1$ contains only two arcs, none of which
connects $F'_a$ with $F''_a$.

After the move of type IC corresponding to  $h=a+1$, we obtain the
situation of Figure 9 in which the Heegaard diagram contains a
pair of complementary handles given by the pair of cycles $E',E''$
and by the cycle $D'_1$, composed by a unique arc connecting $E'$
with $E''$. The deletion of this pair of complementary handles
(Singer move of type III) leads to the genus one Heegaard diagram
drawn in Figure 10, which is the canonical Heegaard diagram of the
lens space $L(2a+1,2r)$.

\vspace{15 pt} {LUIGI GRASSELLI, Department
of Sciences and Methods for Engineering, University of Modena and
Reggio Emilia, 42100 Reggio Emilia, ITALY. E-mail: grasselli.luigi@unimo.it}

\vspace{15 pt} {MICHELE MULAZZANI, Department
of Mathematics, University of Bologna, I-40127 Bologna, ITALY,
and C.I.R.A.M., Bologna, ITALY. E-mail: mulazza@dm.unibo.it}


\bigskip

\begin{figure}[bht]
 \begin{center}
 \includegraphics*[totalheight=10cm]{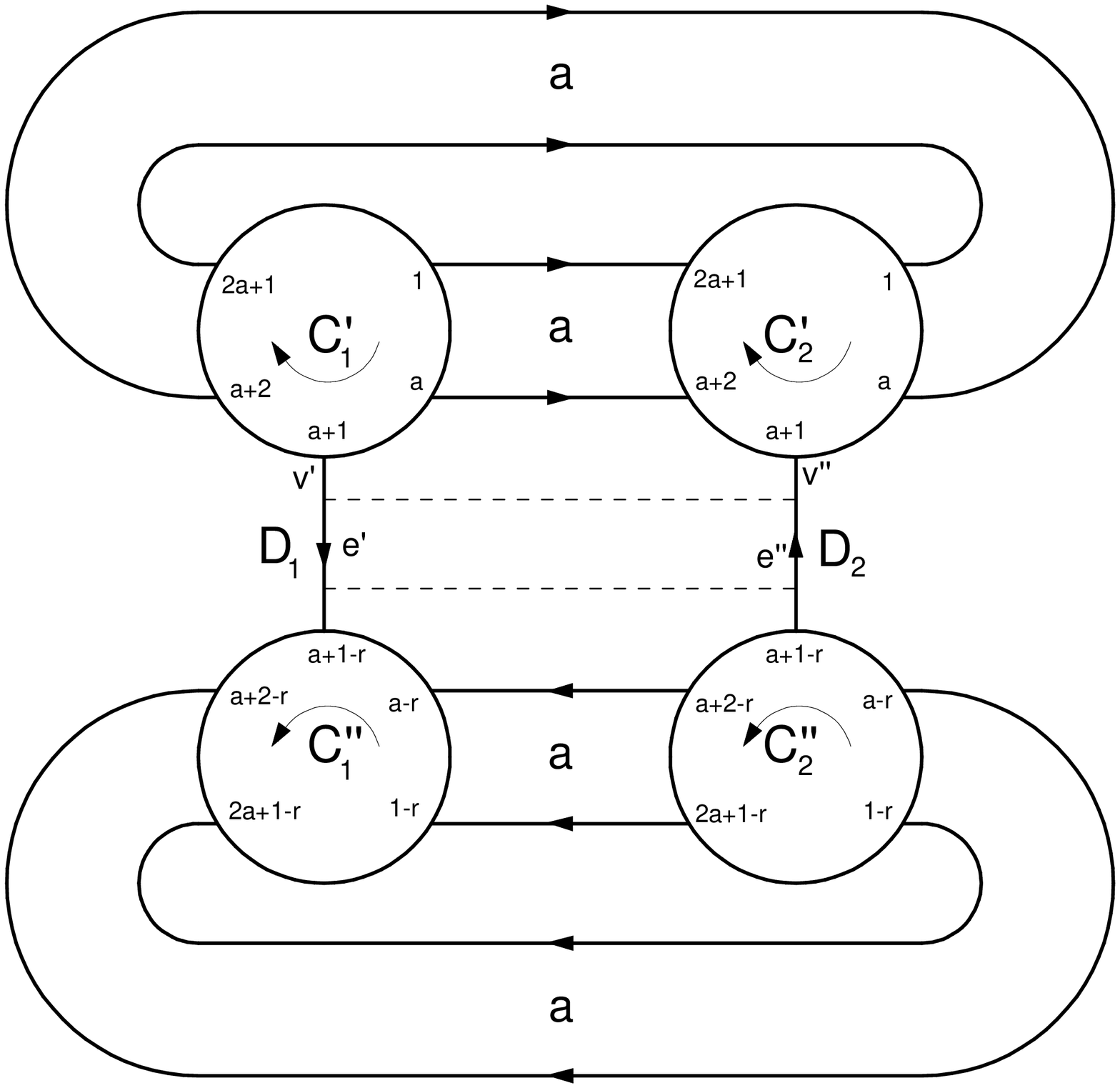}
 \end{center}
 \caption{}

 \label{Fig. 3}

\end{figure}

\bigskip

\begin{figure}[bht]
 \begin{center}
 \includegraphics*[totalheight=11cm]{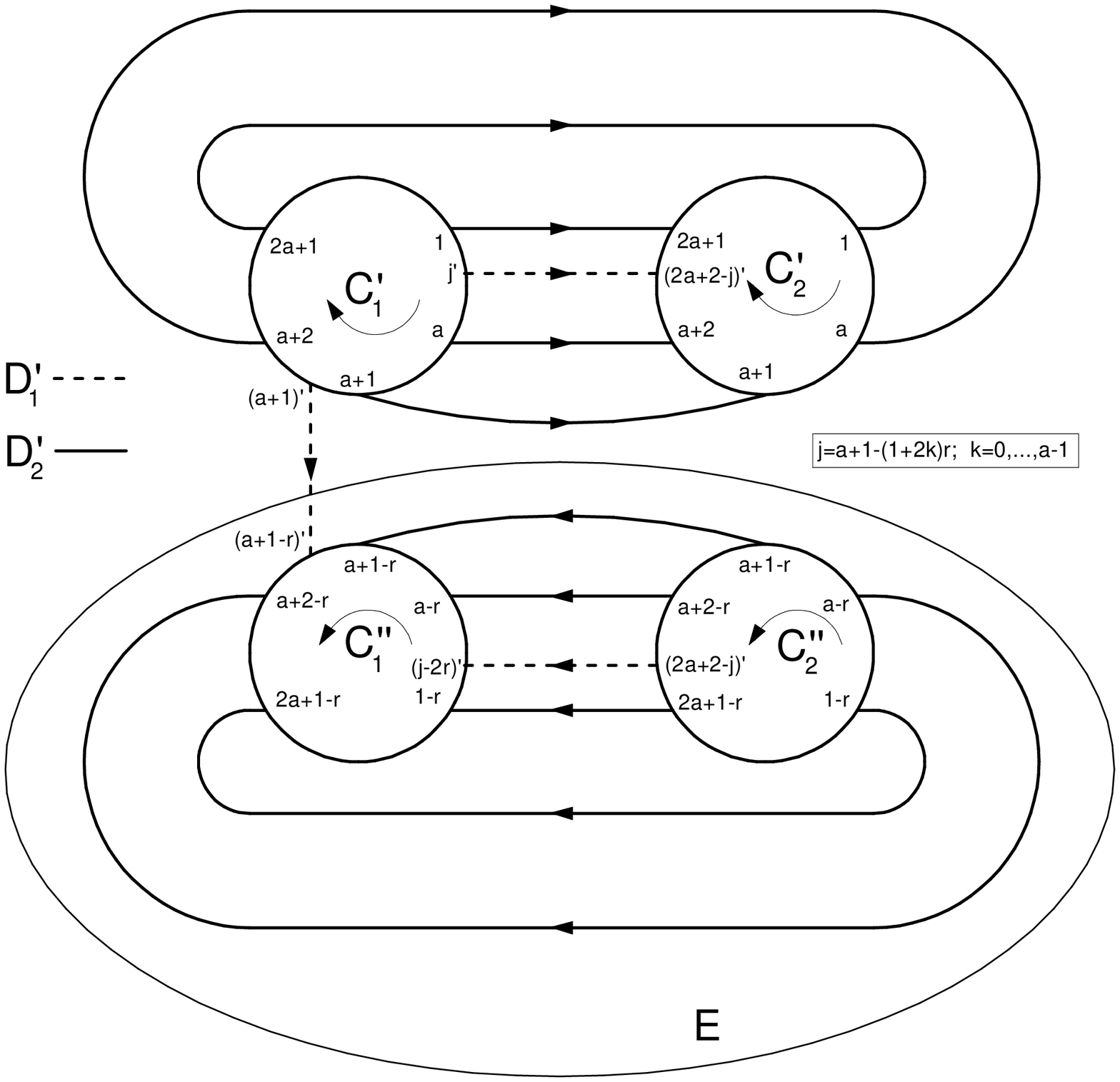}
 \end{center}
 \caption{}

 \label{Fig. 4}

\end{figure}


\bigskip

\begin{figure}[bht]
 \begin{center}
 \includegraphics*[totalheight=7cm]{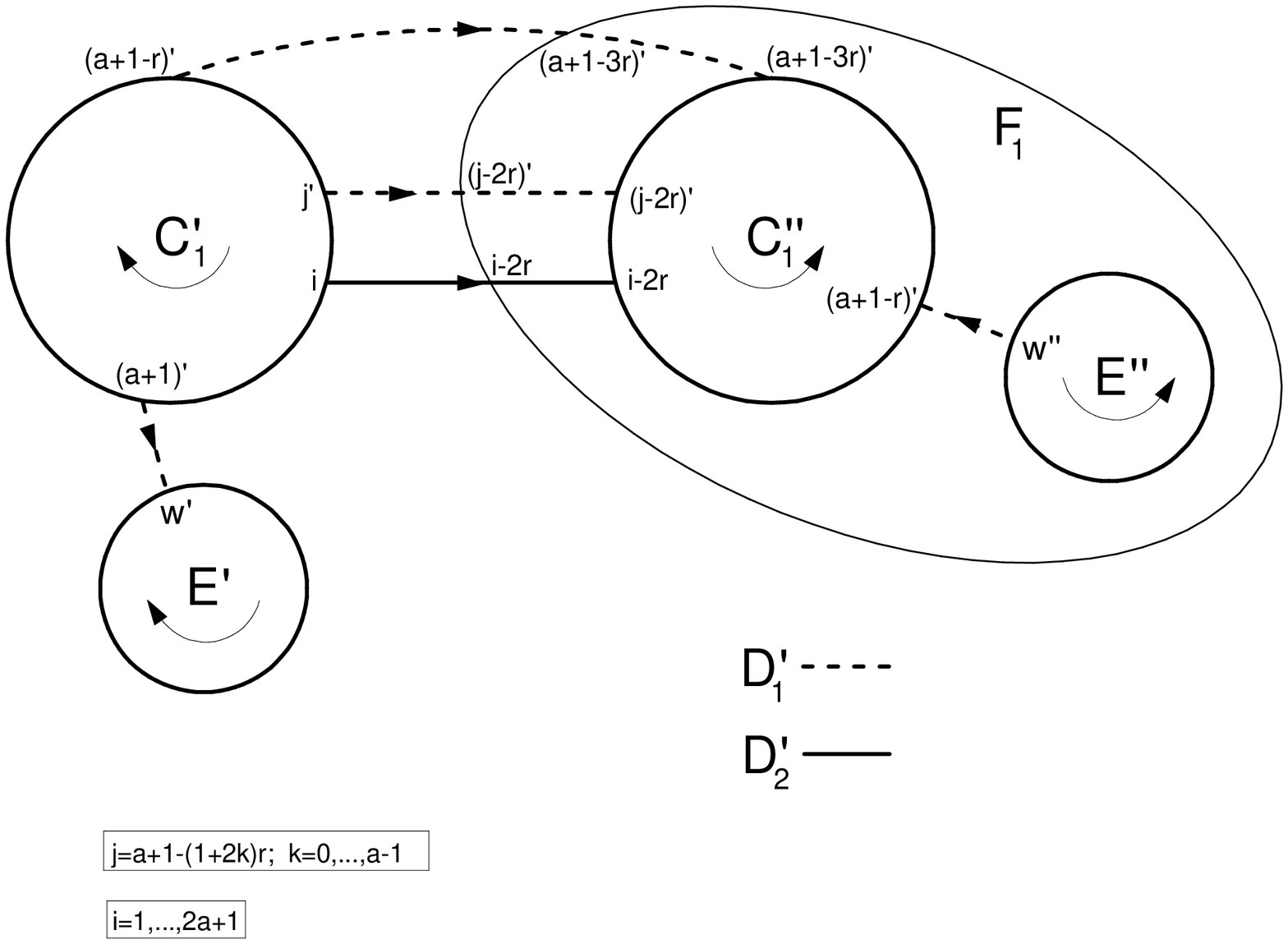}
 \end{center}
 \caption{}

 \label{Fig. 5}

\end{figure}

\bigskip

\begin{figure}[bht]
 \begin{center}
 \includegraphics*[totalheight=9cm]{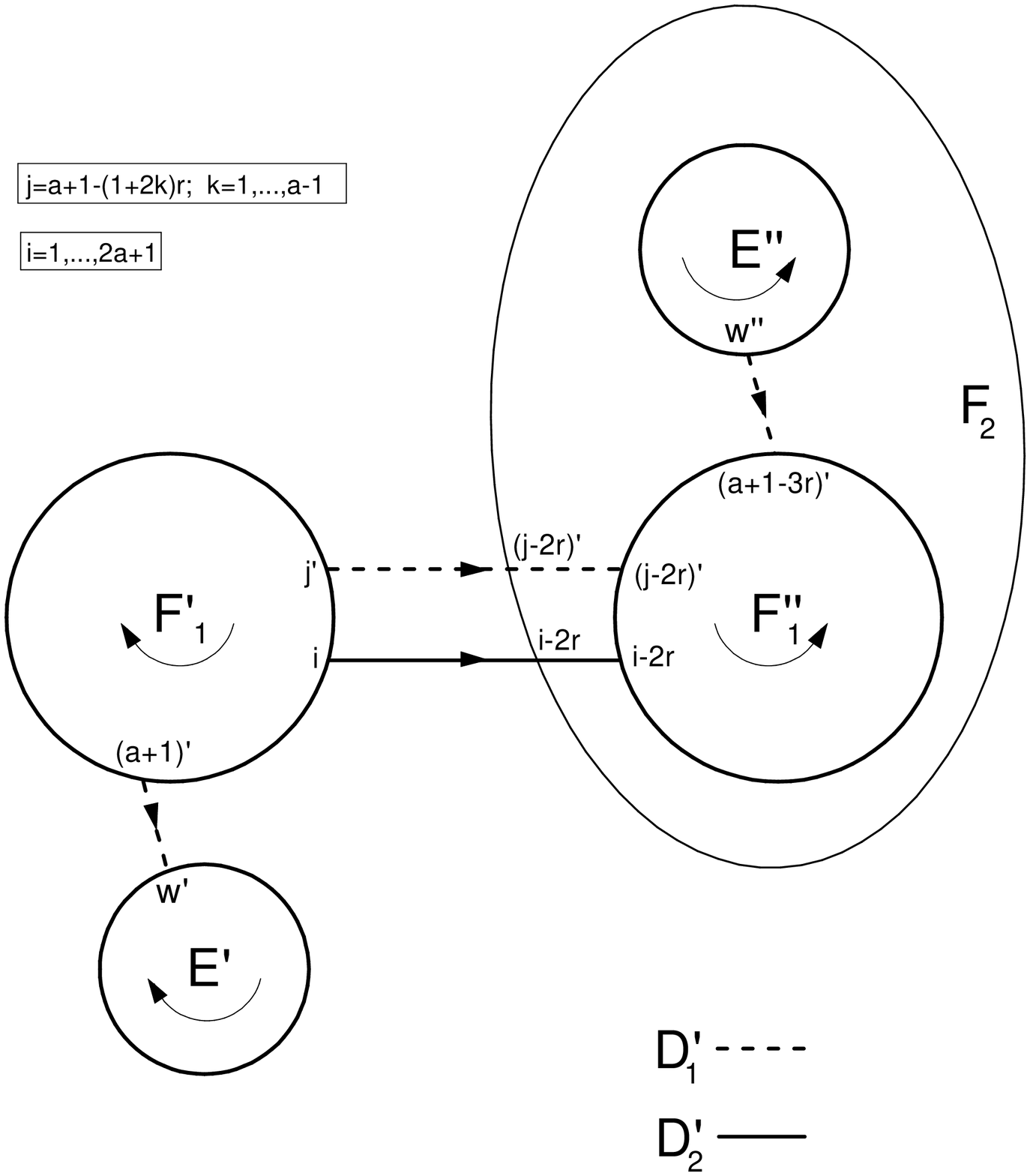}
 \end{center}
 \caption{}

 \label{Fig. 6}

\end{figure}


\bigskip

\begin{figure}[bht]
 \begin{center}
 \includegraphics*[totalheight=7cm]{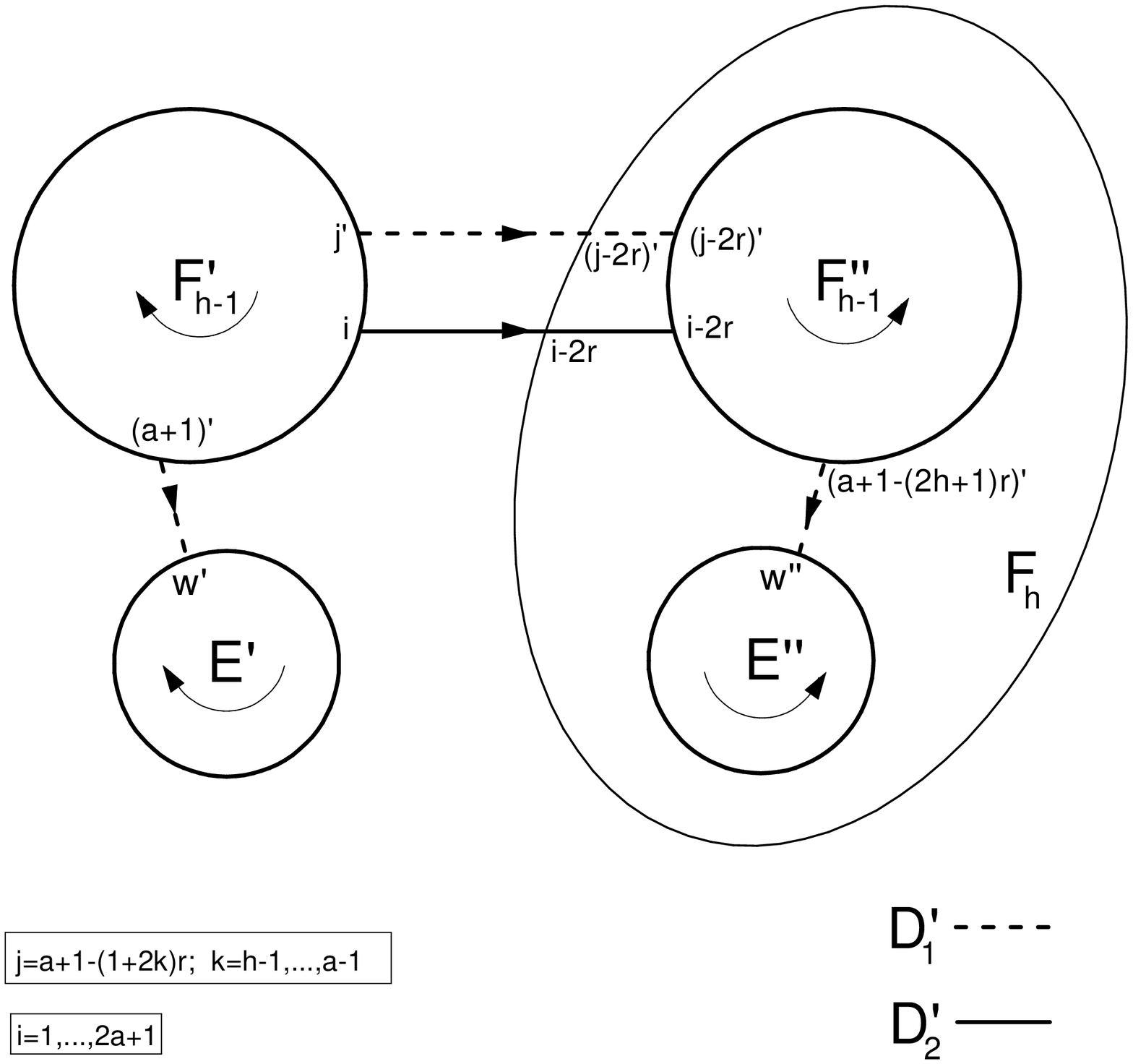}
 \end{center}
 \caption{}

 \label{Fig. 7}

\end{figure}

\bigskip

\begin{figure}[bht]
 \begin{center}
 \includegraphics*[totalheight=8cm]{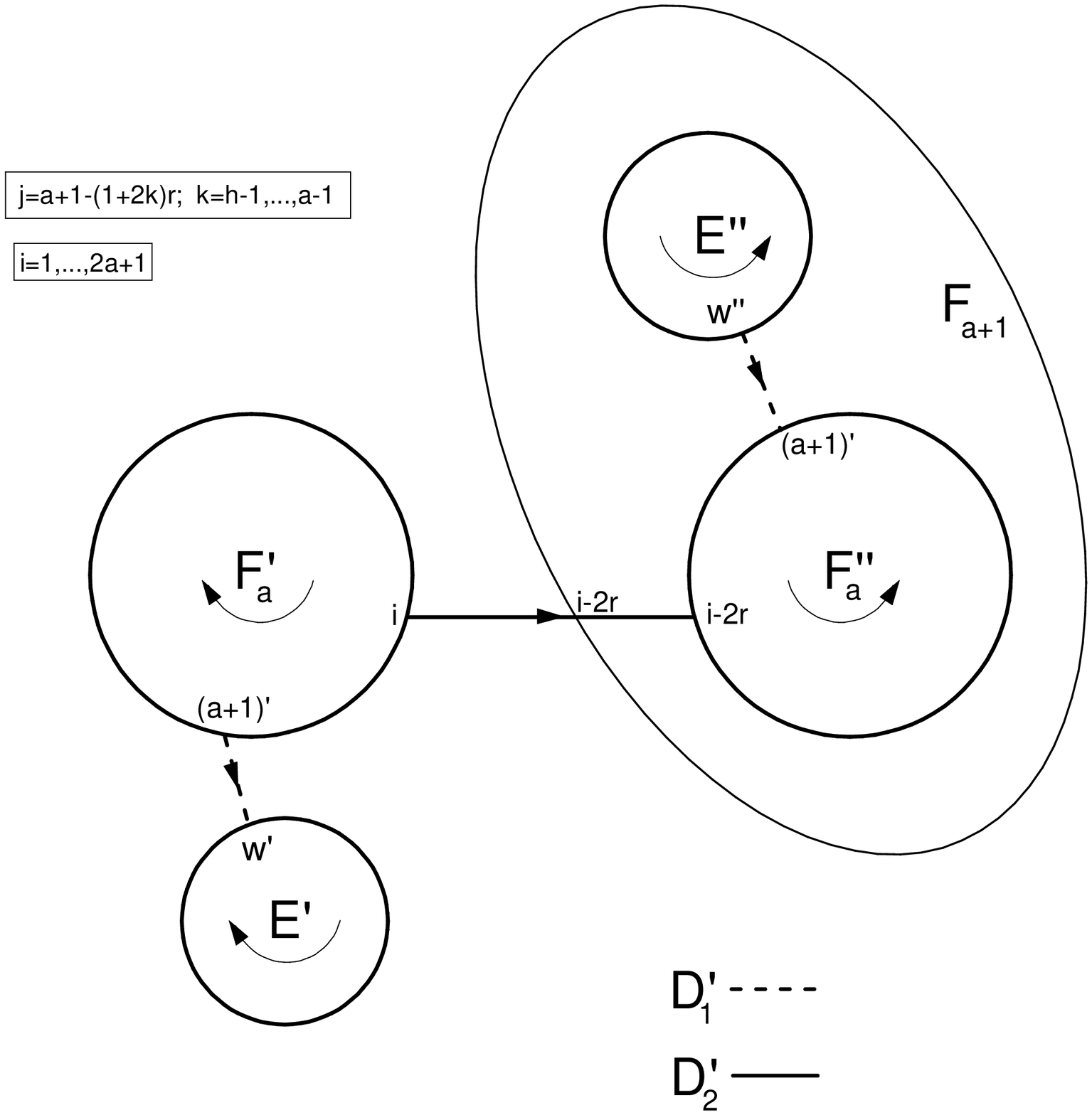}
 \end{center}
 \caption{}

 \label{Fig. 8}

\end{figure}


\bigskip

\begin{figure}[bht]
 \begin{center}
 \includegraphics*[totalheight=7cm]{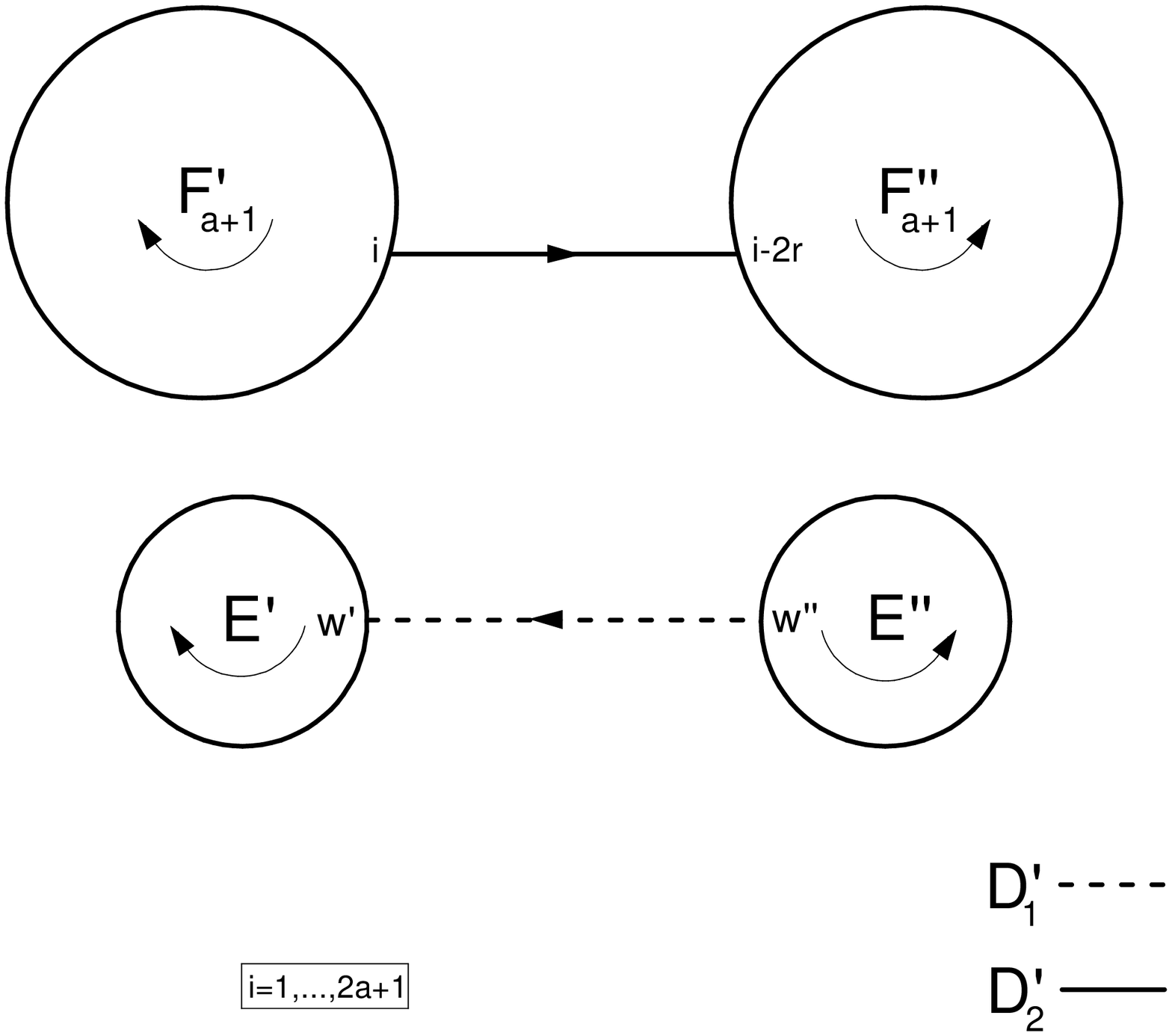}
 \end{center}
 \caption{}

 \label{Fig. 9}

\end{figure}

\bigskip

\begin{figure}[bht]
 \begin{center}
 \includegraphics*[totalheight=3cm]{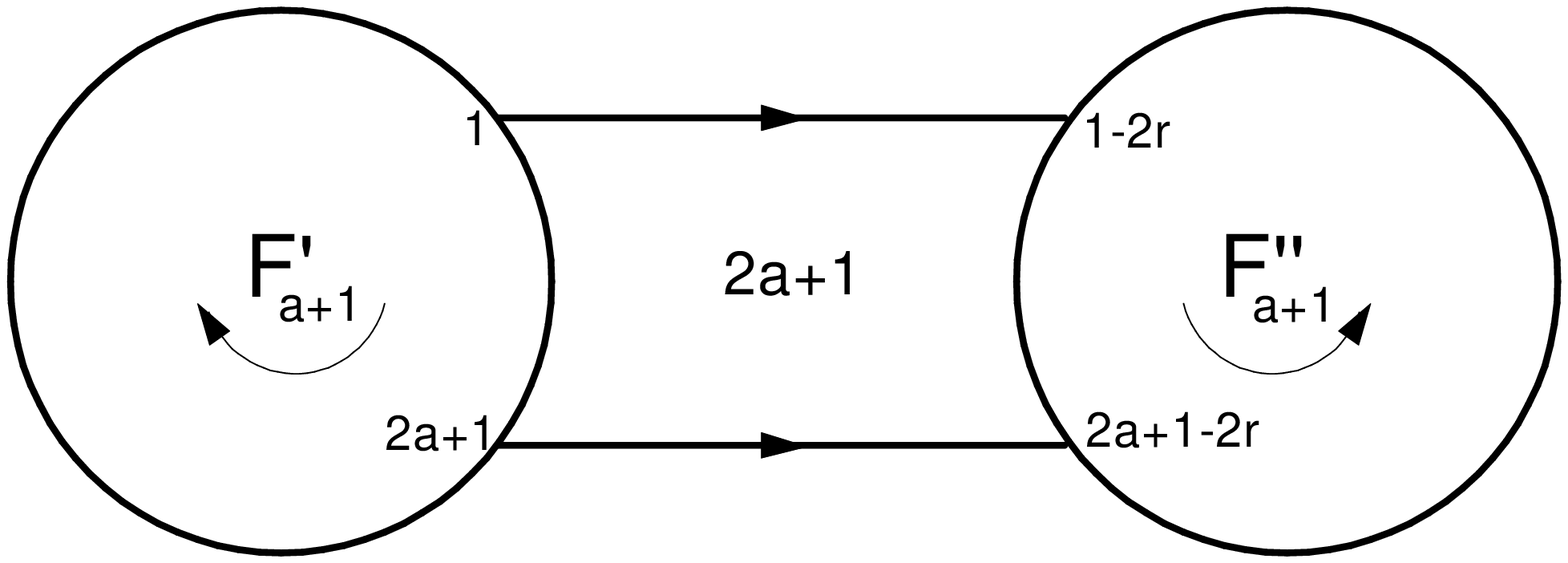}
 \end{center}
 \caption{}

 \label{Fig. 10}

\end{figure}

\end{document}